\documentclass[a4,12pt,epsfig]{elsarticle}
\usepackage{graphicx,epstopdf}
\usepackage{tikz}
\usetikzlibrary{trees}
\usepackage{graphicx}
\usepackage{enumitem}
\usepackage{lipsum}

\usepackage[a4paper, total={5.25in, 8in}]{geometry}

\usepackage{epsfig,epsf,epstopdf}
\usepackage{float}
\usepackage{mathptmx}
\usepackage{amsmath}
\usepackage{amssymb}
\usepackage{amsfonts}
\usepackage{mathtools}
\graphicspath{ {plots/} }
\usepackage{euscript}
\usepackage{textcomp}
\usepackage[subtle]{savetrees}
\usepackage{tikz, pgfplots}
\usetikzlibrary{shapes,arrows}
\usepackage{fancyhdr}
\usetikzlibrary{shadows,calc}
\pgfmathdeclarefunction{gauss}{2}{%
	\pgfmathparse{1/(#2*sqrt(2*pi))*exp(-((x-#1)^2)/(2*#2^2))}%
}

\abovedisplayshortskip
\belowdisplayshortskip
\abovedisplayskip
\belowdisplayskip

\renewcommand{\baselinestretch}{.9}
\renewcommand{\theequation}{\arabic{section}.\arabic{equation}}
\newcommand{\bea}{\begin{eqnarray}}
	\newcommand{\eea}{\end{eqnarray}}
\newcommand{\nbea}{\begin{eqnarray*}}
	\newcommand{\neea}{\end{eqnarray*}}

\newtheorem{theorem}{Theorem}[section]

\newtheorem{lemma}{Lemma}[section]
\newtheorem{example}{Example}[section]
\newtheorem{remark}{Remark}[section]

\begin{document}

\renewcommand{\baselinestretch}{1.25}
\renewcommand{\theequation}{\arabic{section}.\arabic{equation}}
	\bibliographystyle{plain}
	\newcommand{\bfmN}{{\mbox{\boldmath{$N$}}}}
	\newcommand{\bfmx}{{\mbox{\boldmath{$x$}}}}
	\newcommand{\bfmv}{{\mbox{\boldmath{$v$}}}}
	\newcommand{\se}{\setcounter{equation}{0}}
	\newcommand{\vtwo}{\vskip 4ex}
	\newcommand{\vthree}{\vskip 6ex}
	\newcommand{\vfour}{\vspace*{8ex}}
	\newcommand{\hone}{\mbox{\hspace{1em}}}
	\newcommand{\hon}{\mbox{\hspace{1em}}}
	\newcommand{\htwo}{\mbox{\hspace{2em}}}
	\newcommand{\hthree}{\mbox{\hspace{3em}}}
	\newcommand{\hfour}{\mbox{\hspace{4em}}}
	\newcommand{\von}{\vskip 1ex}
	\newcommand{\vone}{\vskip 2ex}
	\newcommand{\n}{\mathfrak{n} }
	\newcommand{\m}{\mathfrak{m} }
	\newcommand{\q}{\mathfrak{q} }
	\newcommand{\aF}{\mathfrak{a} }
	
	\newcommand{\kl}{\mathcal{K}}
	\newcommand{\p}{\mathcal{P}}
	\newcommand{\Lt}{\mathcal{L}}
	\newcommand{\bv}{{\mbox{\boldmath{$v$}}}}
	\newcommand{\bc}{{\mbox{\boldmath{$c$}}}}
	\newcommand{\bx}{{\mbox{\boldmath{$x$}}}}
	\newcommand{\br}{{\mbox{\boldmath{$r$}}}}
	\newcommand{\bs}{{\mbox{\boldmath{$s$}}}}
	\newcommand{\bb}{{\mbox{\boldmath{$b$}}}}
	\newcommand{\ba}{{\mbox{\boldmath{$a$}}}}
	\newcommand{\bn}{{\mbox{\boldmath{$n$}}}}
	\newcommand{\bp}{{\mbox{\boldmath{$p$}}}}
	\newcommand{\by}{{\mbox{\boldmath{$y$}}}}
	\newcommand{\bz}{{\mbox{\boldmath{$z$}}}}
	\newcommand{\be}{{\mbox{\boldmath{$e$}}}}
	
	\newcommand{\bP}{{\mbox{\boldmath{$P$}}}}
	
	\newcommand{\M}{\mathcal{M}}
	\newcommand{\R}{\mathbb{R}}
	\newcommand{\Q}{\mathbb{Q}}
	\newcommand{\Z}{\mathbb{Z}}
	\newcommand{\N}{\mathbb{N}}
	\newcommand{\C}{\mathbb{C}}
	\newcommand{\xar}{\longrightarrow}
	\newcommand{\ov}{\overline}
	\newcommand{\rt}{\rightarrow}
	\newcommand{\om}{\omega}
	\newcommand{\wh}{\widehat }
	\newcommand{\wt}{\widetilde }
	\newcommand{\g}{\Gamma}
	\newcommand{\lm}{\lambda}
	
	\newcommand{\eN}{\EuScript{N}}
	\newcommand{\ncom}{\newcommand}
	\newcommand{\norm}{\|\;\;\|}
	\newcommand{\inp}[2]{\langle{#1},\,{#2} \rangle}
	\newcommand{\nrm}[1]{\parallel {#1} \parallel}
	\newcommand{\nrms}[1]{\parallel {#1} \parallel^2}
	\title{Exploring Chebyshev Polynomial Approximations:\\ Error Estimates for Functions of Bounded Variation}
\author{S. Akansha\\
	{ Department of Mathematics}\\{ Manipal Institute of Technology}\\{Manipal Academy of Higher Education - 576104, India.}}
\date{}
\maketitle{}
\begin{center}{\bf Abstract}\end{center}
Approximation theory plays a central role in numerical analysis, undergoing continuous evolution through a spectrum of methodologies. Notably, Lebesgue, Weierstrass, Fourier, and Chebyshev approximations stand out among these methods. However, each technique possesses inherent limitations, underscoring the critical importance of selecting an appropriate approximation method tailored to specific problem domains. This article delves into the utilization of Chebyshev polynomials at Chebyshev nodes for approximation. For sufficiently smooth functions, the partial sum of Chebyshev series expansion offers optimal polynomial approximation, rendering it a preferred choice in various applications such as digital signal processing and graph filters due to its computational efficiency. In this article, we focus on functions of bounded variation, which find numerous applications across mathematical physics, hyperbolic conservations, and optimization. We present two optimal error estimations associated with Chebyshev polynomial approximations tailored for such functions. To validate our theoretical assertions, we conduct numerical experiments. Additionally, we delineate promising future avenues aligned with this research, particularly within the realms of machine learning and related fields.



\section{Introduction}
Polynomial approximation serves as a fundamental method across various domains of numerical analysis \cite{ame-14a, but-16a}. Not only does it provide a robust tool for approximating complex functions, but it also plays a crucial role in numerical integration and solving differential and integral equations. The Lagrange interpolation polynomial at Chebyshev points of the first or second kind has been observed to mitigate the Runge phenomenon \cite{mas-han_03a}, surpassing interpolants at equally spaced points. Moreover, the accuracy of approximation exhibits rapid enhancement with an increase in the number of interpolation points \cite{tre-13a, mei-12a}. Functions of bounded variation hold significant importance in various branches of mathematical physics, optimization \cite{att_but_mic-14a}, free-discontinuity problems \cite{gar-01a}, and hyperbolic systems of conservation laws \cite{daf-05a}. Additionally, these functions find application in image segmentation and related models \cite{bou_cha-00a}. However, despite their relevance, the theory of numerical approximations for such functions remains relatively underdeveloped, primarily due to the inherent singularities they exhibit.

Considerable research has been devoted to approximating non-smooth functions through the lens of decay estimates of series coefficients. Xiang \cite{xia-18a} delves into the decay behavior of coefficients in polynomial expansions of functions with limited regularity, focusing particularly on Jacobi and Gegenbauer polynomial series. The aim is to establish optimal asymptotic results for the decay of these coefficients. The study investigates how the decay rate varies for functions exhibiting both interior and boundary singularities, considering various parameters. More recently, Wang \cite{wan_23a} tackles the issue of error localization in Chebyshev spectral methods when confronted with functions containing singularities. The study initiates with a pointwise error analysis for Chebyshev projections of such functions, revealing that the convergence rate away from the singularity is faster than at the singularity itself by one power of $\frac{1}{x}$. This analysis provides a rigorous elucidation of the observed error localization phenomenon. The findings suggest that Chebyshev spectral differentiations generally outperform their counterparts except near singularities, where the latter exhibit faster convergence. Notably, when the singularity lies within the interval, Chebyshev spectral differentiations achieve even faster convergence rates in the maximum norm.

Liu et al. \cite{liu_19a} introduce a novel theoretical framework grounded in fractional Sobolev-type spaces, leveraging Riemann-Liouville fractional integrals/derivatives for optimal error estimates of Chebyshev polynomial interpolation for functions with limited regularity. Key components include fractional integration by parts and generalized Gegenbauer functions of fractional degree (GGF-Fs). This framework facilitates the estimation of the optimal decay rate of Chebyshev expansion coefficients for functions with singularities, leading to enhanced error estimates for spectral expansions and related operations. In a separate study, Wang \cite{wan_23b} focuses on deriving error bounds for Legendre approximations of differentiable functions using Legendre coefficients by Hamzehnejad \cite{ham_22a}. Additionally, Xie \cite{xie_23a} recently obtained bounds for Chebyshev polynomials with endpoint singularities. Zhang and Boyd \cite{zha_23a} derived estimates for weak endpoint singularities, while Zhang \cite{zha_21a, zha_23b} obtained bounds for logarithmic endpoint singularities. In \cite{occ-the_21a}, the focus lies on a specialized filtered approximation technique that generates interpolation polynomials at Chebyshev zeros using de la Vallée Poussin filters. The aim is to approximate locally continuous functions equipped with weighted uniform norms. Ensuring that the associated Lebesgue constants remain uniformly bounded is crucial for this endeavor.

The methodologies discussed above primarily concentrate on Chebyshev interpolation \cite{occ-ram-the_22a, occ-the_21a, zha_21a, de-mar-occ_21a}, yielding results with exact Chebyshev series coefficients. However, computing exact series coefficients poses a general challenge and proves impractical for numerical algorithms, diminishing their utility in practical applications. Furthermore, these approaches usually involve Jacobi, Gegenbauer, and Legendre polynomials on fixed intervals where the respective series' basis functions are defined. Such limitations highlight the need for more versatile and efficient approximation methods in numerical analysis.
Addressing this gap necessitates the utilization of efficient approximation techniques. Chebyshev polynomials, renowned for their versatility and effectiveness across diverse fields such as digital signal processing \cite{shu_van_fro-11a}, spectral graph neural networks \cite{lee-jen_98a, he-wei_22a, li-wan_24a}, image processing \cite{occ-ram-the_22a} and graph signal filtering \cite{sta_man_dak-19a, tse-lee_21a}, present a promising avenue for approximating functions of bounded variation.

Truncated Chebyshev expansions have proven capable of yielding minimax polynomial approximations for analytic functions \cite{ged-78a}. Our objective is to employ these polynomials not only for approximating functions of bounded variation but also for conducting a comprehensive convergence analysis of Chebyshev polynomial approximation techniques. At the core of our convergence analysis lies the estimation of Chebyshev coefficients decay. We leverage two recently established decay estimates: Majidian's decay estimate for Chebyshev series coefficients of functions defined on the interval $[-1,1]$, subject to specific regularity conditions \cite{maj_17a}; and a sharper decay estimate demonstrated by Xiang \cite{xia_18a}, with a more relaxed smoothness assumption on the function. In our pursuit of convergence results, we take an initial step by extending these decay bounds to encompass Chebyshev series coefficients of functions defined on the interval $[a,b]$.
\section{Preliminaries}
\subsection{Chebyshev Series Expansion}\label{prelsec:chebyshevseries}
The Chebyshev polynomial of the first kind, denoted as $T_j(t)$ for a given integer $j \geq 0$, is defined as:
\begin{equation}\label{eq:chePoly}
	T_j(t) = \cos(j\theta),
\end{equation}
where $\theta = \cos^{-1}(t)$ and $t \in [-1,1]$. Notably, $T_j(t)$ is a polynomial of degree $j$ in the variable $t$. These polynomials exhibit orthogonality with respect to the weight function $\omega(t) = \frac{1}{\sqrt{(1-t^2)}},$ within the interval $[-1,1]$. Specifically, they satisfy the following orthogonality relations:
\begin{equation*}
	\int\limits_{-1}^{1} \omega(s)T_p(s) T_q(s)ds = \begin{cases}
		0 &\quad\text{if } p \neq q\\
		\pi &\quad\text{if } p = q = 0\\
		\frac{\pi}{2} &\quad\text{if } p=q\neq 0.
	\end{cases}
\end{equation*}
The Chebyshev series expansion of a function $f:[-1,1]\rightarrow \mathbb{R}$ is expressed as follows:
\begin{equation}\label{eq:CheExp}
	f(t) = \frac{c_0}{2} + \sum_{j=1}^{\infty} c_jT_j(t), \quad \text{where} \quad c_j = \dfrac{\langle\,f\,,\,T_j\,\rangle_\omega}{\|T_j\|^2_{\omega}},
\end{equation}
and 
\begin{equation*}
	\langle\,f\,,\,T_j\,\rangle_\omega = \int\limits_{-1}^{1} \omega(s) f(s) T_j(s) ds.   
\end{equation*}
The norm $\|T_j\|_\omega$ is computed as:
\begin{equation}\label{eq:Orthogonality}
	\|T_j\|_\omega = \big<\,T_j, T_j\,\big>_\omega^{\frac{1}{2}} =  \begin{cases} 
		\sqrt{\pi} & j = 0 \\
		\sqrt{\pi /2} & j \neq 0.
	\end{cases}
\end{equation}
Hence, the Chebyshev coefficients $c_j$ can be obtained using the integral form:
\begin{equation}\label{eq:Checoff}
	c_j = \frac{2}{\pi}\int\limits_{-1}^{1}f(s)T_j(s)\omega(s) ds.
\end{equation} 
Given the difficulty in accurately evaluating the integral \eqref{eq:Checoff} for general functions, we resort to employing the Gauss Chebyshev quadrature rule to approximate $c_j$, the $j$th coefficient of the series.

\subsubsection{Gauss Chebyshev quadrature formula}\label{prelsec:gaussquad}
Quadrature methods are renowned for numerically computing definite integrals of the type presented in \eqref{eq:Checoff}. The Gauss-Chebyshev quadrature formula, a variant of Gaussian quadrature employing the weight function $\omega$ and $n$ Chebyshev points, provides an explicit formula for numerical integration (see \cite{mas-han_03a,tam_lop-07a,riv_74a}):
\begin{equation}\label{eq:quadformula}
	\int\limits_{-1}^{1} \omega(s)F(s)ds \approx \frac{\pi}{n}\sum_{l=1}^{n}F(t_l),   
\end{equation}
where $t_1,t_2,\ldots,t_n$ represent the $n$ roots of a Chebyshev polynomial $T_n(t)$ of degree $n$, also known as Chebyshev points, defined as:
\begin{equation}\label{eq:ChePts}
	t_l = \cos\left(\dfrac{\left(2l-1\right) \pi}{2n}\right), \quad l=1, 2, \ldots,n. 
\end{equation}
Leveraging the quadrature formula \eqref{eq:quadformula}, we can readily approximate Chebyshev series coefficients \eqref{eq:Checoff} for any function using the formula provided by Rivlin \cite{riv_74a} (page 148):
\begin{equation}\label{eq:approxcoeff}
	c_k \approx \frac{2}{n}\sum_{l=1}^{n}f(t_l)T_k(t_l) =:c_{k,n}.
\end{equation}
Here, $c_{k,n}$ denotes the approximated Chebyshev coefficients using $n$ quadrature points.

\subsection{Notations}
We denote the Chebyshev series expansion of a function $f\in L^2_\omega[a,b]$ by $\mathsf{C}_\infty[f](x)$, defined as:
$$\mathsf{C}_\infty[f](x) := \sum_{j=0}^{\infty} c_j T_j(\mathsf{G}^{-1}(x)),$$
where $\mathsf{G}:[-1,1]\rightarrow [a,b]$ is a bijection map given by:
$$\mathsf{G}(t) = a+\frac{(b-a)}{2}(t+1),~t\in [-1,1].$$
The Chebyshev coefficients $c_j$ are calculated as:
$$c_j =  \frac{2}{\pi}\int_{-1}^1 
\frac{f(\mathsf{G}(t)) T_j(t)}{\sqrt{1-t^2}}dt .$$
Utilizing the change of variable $t = \cos \theta$, we express $c_j$ as:
\begin{equation}\label{eq:c_jcostheta}
	c_j = \dfrac{2}{\pi}\int_{0}^{\pi} f(\mathsf{G}(\cos \theta)) \cos j\theta d\theta.
\end{equation}
The $d^{\rm th}$ partial sum, $C_d[f](x)$, approximates the function $f$ at a point $x\in [a,b]$, given by:
\bea\label{nthPartialSum.eq}
C_d[f](x) := \sideset{}{'}\sum_{j=0}^{d} c_j T_j(\mathsf{G}^{-1}(x))
\eea
In our results, we use:
\bea\label{CCoeff.Approx.eq}
C_{d,n}[f](x) := \sideset{}{'}\sum_{j=0}^{d} c_{j,n} T_j(\mathsf{G}^{-1}(x))
\eea
to denote the corresponding approximation of $f$ using $n$ quadrature points.

Additionally, we represent the Chebyshev series expansion of $f$, approximated with $n$ quadrature points, as:
\begin{equation}
	\mathsf{C}_{\infty,n}[f](t) := \sideset{}{'}\sum_{j = 0}^{\infty} c_{j,n}T_j(t),
\end{equation}
The following lemmas are used in deriving the required error estimates.
\begin{lemma}\label{Ident.ChebApproxAccur.lemma}
	For a given positive integer $n$, we have
	$$c_{k,n} - c_k = \sum_{j=1}^\infty (-1)^j\big(c_{2jn-k} + c_{2jn+k}\big),$$
	for any integer $k$ such that $0\le k< 2n$.
\end{lemma}
\begin{lemma}\label{Estimate.ChebCoeff.lem}
	For $0\le d<2n$, we have
	$$\|C_d[f] - C_{d,n}[f]\|_1\le (b-a)\sum_{j=1}^\infty \sum_{i=2jn-d}^{2jn+d}|c_i|,$$
	for any $f\in L^1[a,b].$
\end{lemma}
Proof of these lemmas are given in Appendix-A.
\section{Decay bounds for Chebyshev Coefficients}\label{decay.bounds.sec}
In this section, we extend the decay bounds established in prior works by Majidian \cite{maj_17a} and Xiang \cite{xia_18a}. This generalization is pivotal for numerous applications. In practical scenarios, the function to be approximated may not always reside solely within the domain $[-1,1]$. Moreover, in various applications, local schemes \cite{car-hoc-stu_20a} or piecewise approximation \cite{glu-sie-tho_24a, van_14a} are preferred over global ones. In such cases, decay estimates on a general domain become imperative.
\begin{theorem}\label{thm:coeffbound}
	For some integer $k\ge 0$, let $f$, $f'$, $\ldots,$ $f^{(k-1)}$ be absolutely continuous on the interval $[a, b]$. 
	If $V_k:=\|f^{(k)}\|_T < \infty,$ where
	\begin{equation}\label{eq:chebnormscale}
		\|f\|_T := \int_{0}^{\pi} \left|{f'\left(\mathsf{G}(\cos \theta )\right)}\right|d\theta,
	\end{equation}
	then for $j\ge k+1$ and for some $s\ge 0$, we have
	\begin{equation}\label{eq:coeffbound1}
		|c_j| \leq 
		\left\{
		\begin{array}{@{}l@{\thinspace}l}
			\left(\dfrac{b-a}{2}\right)^{2s+1}  \dfrac{2V_k}{\pi\displaystyle \prod_{i = -s}^{s}(j + 2i)},&
			\mbox{ if }~ k = 2s,\medskip\\
			\left(\dfrac{b-a}{2}\right)^{2s+2}  \dfrac{2V_k}{\pi\displaystyle \prod_{i = -s}^{s+1}(j + 2i-1)},&
			\mbox{ if }~ k = 2s+1,
		\end{array}
		\right. 
	\end{equation}
	where $c_k$, $k=0,1,\ldots$, are the Chebyshev coefficients of $f$.
\end{theorem}
\textbf{Proof: } For a given nonnegative integer $r$, we define:
\begin{equation}\label{eq:c_j^r}
	c_j^{(r)} = \dfrac{2}{\pi}\int_{0}^{\pi} f^{(r)}(\mathsf{G}(\cos \theta)) \cos j\theta d\theta,
\end{equation}
with the understanding that $c_j^{(0)} = c_j$. In dealing with non-smooth functions, we must utilize the weak derivative (distributional derivative) on the right-hand side of the above expression, if it exists.
Employing integration by parts in \eqref{eq:c_j^r}, we can express $c_j^{(r)}$ as:
\begin{equation}\label{eq:c_j^ridentity}
	c_j^{(r)} = \dfrac{(b-a)}{4j}\left(c_{j-1}^{(r+1)} - c_{j+1}^{(r+1)}\right),
\end{equation}
for $j = 1, 2, \ldots$. In order to prove the required estimate, we prove the following general inequality 
\begin{equation}\label{eq:coeffboundgeneral}
	|c_j^{(k-m)}| \leq \dfrac{2V_k}{\pi}\left\{
	\begin{array}{@{}l@{\thinspace}l}
		\left(\dfrac{b-a}{2}\right)^{m+1} \dfrac{1}{\displaystyle \prod_{i = -s}^{s}(j+2i)},  & \text{ if } m = 2s,s\ge0\\
		
		\left(\dfrac{b-a}{2}\right)^{m+1}  \dfrac{1}{\displaystyle \prod_{i = -s}^{s+1}(j+2i-1)}, & \text{ if } m = 2s+1, s\ge 0 \\
	\end{array}
	\right. 
\end{equation}
for $m = 0,....,k$ and $j \geq m+1$ using induction on $m$. Then $m=k$ gives our required result. 
First let us claim that \eqref{eq:coeffboundgeneral} holds for $m = 0$. From \eqref{eq:c_j^ridentity} and \eqref{eq:c_j^r}, we have
\begin{equation*}
	|c_j^{(k)}|  \leq \dfrac{b-a}{4j}\left(|c_{j-1}^{(k+1)}| + |c_{j+1}^{(k+1)}|\right)\le \dfrac{(b-a)}{j\pi} V_k.
\end{equation*}
This is precisely the inequality \eqref{eq:coeffboundgeneral} for $m=0.$
Let us assume that the inequality \eqref{eq:coeffboundgeneral} holds for $m = 2s$ for some $s\ge 1$.  Then for $m=2s+1$ (odd), we have
\begin{align*}
	|c_j^{(k-2s-1)}| & \leq \dfrac{b-a}{4j}\left(|c_{j-1}^{(k-2s)}| + |c_{j+1}^{(k-2s)}|\right) 
\end{align*}
Using the assumption that the inequality \eqref{eq:coeffboundgeneral} holds for $m=2s$, we can write
\begin{align*}
	|c_j^{(k-2s-1)}| & \leq \dfrac{b-a}{4j}\left(\left(\dfrac{b-a}{2}\right)^{2s+1} \frac{2V_k}{\pi}\left[\dfrac{1}{\displaystyle \prod_{i = -s}^{s}(j+2i-1)}  +
	\dfrac{1}{\displaystyle \prod_{i = -s}^{s}(j+2i+1)}\right]
	\right).
\end{align*}
Up on simplifying the right hand side, we get
\begin{align*}
	|c_j^{(k-2s-1)}|&\le
	\left(\dfrac{b-a}{2}\right)^{2s+2}\frac{2V_k}{\pi\displaystyle \prod_{i = -s}^{s+1}(j+2i-1)},
\end{align*}
which is precisely the required inequality \eqref{eq:coeffboundgeneral} for $m=2s+1.$ Similarly, asssume that the inequality \eqref{eq:coeffboundgeneral} holds for $m = 2s+1$. Then for $m=2s+2$ (even), we have
\begin{align*}
	|c_j^{(k-2s-1)}|&\le
	\left(\dfrac{b-a}{2}\right)^{2s+3}\frac{2V_k}{\pi\displaystyle \prod_{i = -(s+1)}^{s+1}(j+2i)}.
\end{align*}
which is precisely the required inequality \eqref{eq:coeffboundgeneral} for $m=2s+2.$ The proof now follows by induction.
\begin{remark}
	Note that if $f^{(k)}$ is absolutely continuous, then $V_k$ is precisely the total variation of $f^{(k)}$ and hence in this case, the assumption that $V_k$ being finite implies $f^{(k)}$ is of BV on $[a,b]$.  If $f^{(k)}$ involves a jump discontinuity, then one has to necessarily use distribution derivative of $f^{(k)}$ in computing $V_k$.  
\end{remark}
The following lemma is the generalization of a result in \cite{xia-18a}. 
\begin{lemma}\label{thm:c_jformula}
	Let $f$ be a function defined on an interval $[a,b]$ such that for some integer $k\geq 1$, $c_j^{(k)}$ is well defined and $f^{(k)}$ is of bounded variation on $[a,b]$. Then we have
\begin{equation}\label{eq:c_jformula}	
	c_j = \left(\dfrac{b-a}{4}\right)^p \sum_{i=0}^{p} {p \choose i}\dfrac{(-1)^i(j+2i-p)}{(j+i)(j+i-1)\ldots (j+i-p)}c^{(p)}_{(j+2i-p)}, 
\end{equation}
where $j = p,p+1,\ldots$ and $p = 1,2, \ldots,k$. 
\end{lemma}
\begin{theorem}\label{thm:decayc_jinV_k}
	Let $f$ be a function defined on $[a,b]$ such that for some nonnegative integer $k$, $f^{(k)}$ is of bounded variation with $V_k = Var(f^{(k)})<\infty$. Then we have
	\begin{align}
		|c_j^{(k)}| & \leq \dfrac{2V_k}{j\pi}, ~~~~~~~~~~~~~~~~~~~~ j = 1,2,\ldots, ~~~~~~~~~~~~~~~~\mbox{ and }\\
		|c_j| &\leq \dfrac{2V_k}{\pi}\left(\dfrac{b-a}{4}\right)^k\sum_{i=0}^{k} {k \choose i}\dfrac{1}{(j+i)(j+i-1)\ldots (j+i-k)},\label{cjXiang.eq}
	\end{align}
	for $j = k+1,k+2,\ldots$.
\end{theorem}
\textbf{Proof: } Since $f^{(k)}$ is of bounded variation we can write (see Lang \cite{lan-12a})
\begin{equation}\label{Varfk.eq}
	Var\left(f^{(k)}\right) = Var(g_1)+Var(g_2),
\end{equation}
where $g_1$ and $g_2$ are monotonically increasing functions on $[a,b]$. Define $u(\theta):=g_i(\mathsf{G}(\cos\theta))$ which is monotonically decreasing for  $\theta \in [0,\pi]$. Further, we have
\begin{align}\label{eq:definec_n^(k)ing}
	\int_{0}^{\pi}g_i(\mathsf{G}(\cos\theta))\cos j\theta \, d\theta 
	& = -\dfrac{2}{\pi}\int_{0}^{\pi}v(\theta)\cos j\theta \, d\theta. \nonumber
\end{align}
where  $v(\theta) = -u(\theta)$. By second mean value theorem of integral calculus (Apostol \cite{apo-74a}, Theorem 7.37), there exist $x_0 \in [0,\pi]$ such that
\begin{equation}\label{eq:c_n^(k)thm}
	\int_{0}^{\pi}g_i(\mathsf{G}(\cos\theta))\cos j\theta \, d\theta = -\dfrac{2}{\pi}\left(v(0)\int_{0}^{x_0}\cos j\theta d\theta +v(\pi)\int_{x_0}^{\pi}\cos j\theta \, d\theta \right).
\end{equation}
By definition of $v$, we have 
\begin{align*}
	-v(0) &= u(0) = g_i(\mathsf{G}(\cos 0)) = g_i(\mathsf{G}(1)) = g_i(b)\\
	-v(\pi) &= u(\pi) = g_i(\mathsf{G}(\cos \pi)) = g_i(\mathsf{G}(-1)) = g_i(a).
\end{align*}
Substituting in \eqref{eq:c_n^(k)thm} and then integration yields 
\begin{align*}
	\int_{0}^{\pi}g_i(\mathsf{G}(\cos\theta))\cos j\theta d\theta
	& = \dfrac{2}{\pi}\dfrac{g_i(b)-g_i(a)}{j}\sin jx_0.
\end{align*}
Using \eqref{eq:c_j^r} and \eqref{Varfk.eq}, we get \(|c_j^{(k)}| \leq  \dfrac{2V_k}{j\pi},\) for \(j = 1,2,\ldots.\) Consider \eqref{eq:c_jformula} with $p = k$ which gives
\begin{align*}
	c_j &= \left(\dfrac{b-a}{4}\right)^k \sum_{i=0}^{k} {k \choose i}\dfrac{(-1)^i(j+2i-k)}{(j+i)(j+i-1)\ldots (j+i-k)}c^{(k)}_{(j+2i-k)}.
\end{align*}
Taking modulus on both sides and using the above inequality, we get 
\begin{align*}
	|c_j| 
	& \leq \left(\dfrac{b-a}{4}\right)^k \sum_{i=0}^{k} {k \choose i}\dfrac{(j+2i-k)}{(j+i)(j+i-1)\ldots (j+i-k)}\left(\dfrac{2V_k}{(j+2i-k)\pi}\right),
\end{align*}
for $j = k+1,k+2,\ldots,$ which leads to the desired result. Note that, in this case $j = k$ is not defined when $i$ is $0$.

For the well-known decay estimate of the Chebyshev coefficients of a real analytic function we refer to Rivlin \cite{riv_74a} (also see Xiang {\it et al.} \cite{xia-etal_10a}). 

\section{Error Estimate for Chebyshev Approximation}\label{Cheb.Approx.sec}
In this section, we derive $L^1$-error estimates for the Chebyshev approximation of \(f\), utilizing two decay estimates provided in Theorems \ref{thm:coeffbound} and \ref{thm:decayc_jinV_k}. Specifically, we establish the error estimate for the truncated Chebyshev series approximation, relying on the decay estimate \eqref{eq:coeffbound1} of the Chebyshev coefficients, as presented in the following theorem.
\begin{theorem}\label{thm:errestimate}
	Assume the hypotheses of Theorem \ref{thm:coeffbound}. Then for any given integers $n$ and $d$ such that $n-1\ge k\ge 1$ and $k\le d \le 2n-k-1$, we have
	$$\|f - \mathsf{C}_{d,n}[f]\|_1
	\leq \mathbb{T}_{d,n}$$
	where,
	\begin{enumerate}
		\item if $d=n-l$, for some $l =1,2,\ldots, n-k$, then we have
		\begin{equation}\label{eq:totalL1err}
			\mathbb{T}_{d,n}:= \left\{
			\begin{array}{@{}l@{\thinspace}l}
				\left(\dfrac{b-a}{2}\right)^{k+2} \dfrac{4V_k}{k \pi } \left(\Pi_{1,1}(n-l) + \Pi_{1,2}(n-l)\right),  & \text{ if } k = 2s,\medskip\\		
				\left(\dfrac{b-a}{2}\right)^{k+2} \dfrac{4V_k}{k\pi} \left(\Pi_{0,0}(n-l) + \Pi_{0,1}(n-l)\right), & \text{ if } k = 2s+1, \\
			\end{array}
			\right.
		\end{equation}
		\item else if $d = n+l$, for some $l =0,1,\ldots, n-k-1$, then we have
		\begin{equation}
			\mathbb{T}_{d,n} := \left\{
			\begin{array}{@{}l@{\thinspace}l}
				\left(\dfrac{b-a}{2}\right)^{k+2} \dfrac{6V_k}{\pi k } \left(\Pi_{1,0}(n-l) + \Pi_{1,1}(n-l)\right),  & \text{ if } k = 2s,\medskip\\		
				\left(\dfrac{b-a}{2}\right)^{k+2} \dfrac{6V_k}{k\pi} \left(\Pi_{0,-1}(n-l) + \Pi_{0,0}(n-l) \right), & \text{ if } k = 2s+1, \\
			\end{array}
			\right.
		\end{equation}
		for some integer $s\ge 0$, where
		\begin{equation}\label{PI.def}
			\Pi_{\alpha,\beta}(\eta) := \dfrac{1}{\displaystyle \prod_{i = -s}^{s-\alpha}(\eta+2i+\beta)},~~\alpha=0,1,~~\beta=-1,0,1,2.
		\end{equation}
	\end{enumerate}
\end{theorem}
\textbf{Proof:} 
We have 
$$
\big\|f-C_{d,n}[f]\big\|_1
\le
\big\|f-C_{d}[f]\big\|_1 + \big\|C_{d}[f]-C_{d,n}[f]\big\|_1.
$$
For estimating the second term on the right hand side of the above inequality, we use  the well-known result
(see Fox and Parker \cite{fox-par_68a}) 
$$c_{d,n} - c_d = \sum_{j=1}^\infty (-1)^j\big(c_{2jn-d} + c_{2jn+d}\big),$$
for $0\le d< 2n$.  Using this property, with an obvious rearrangements of the terms in the series, we can obtain
\begin{align}\label{Est01.eq}
	\|f - \mathsf{C}_{d,n}[f]\|_1
	& \leq (b-a)\mathcal{E},
\end{align}
where
$$\mathcal{E}:=\left\{\displaystyle{\sum_{j=d+1}^\infty} |c_j|+\sum_{j=1}^\infty \sum_{i=2jn-d}^{2jn+d}|c_i|\right\}.$$
By adding some appropriate positive terms, we can see that (also see Xiang {\it et al.} \cite{xia-etal_10a})
\begin{equation}\label{E.estimate}
	\mathcal{E} \le 
	\left\{\begin{array}{ll}
		2\displaystyle{\sum_{i=n-l+1}^{\infty}}|c_i|,&\mbox{for } d=n-l,~l =1,2,\ldots, n,\\
		3\displaystyle{\sum_{i=n-l}^{\infty}}|c_i|,&\mbox{for } d=n+l,~l=0,1,\ldots,n-1.
	\end{array}\right.
\end{equation}
We restrict the integer $d$ to $k\le d \le 2n-k-1$ so that the decay estimate in Theorem \ref{thm:coeffbound} can be used.   Now using the telescopic property of the resulting series (see also Majidian \cite{maj_17a}), we can arrive at the required estimate.
\begin{remark}\label{UB.behavior.remark}
	From the above theorem, we see that for a fixed $n$ (as in the hypotheses), the upper bound $\mathbb{T}_{d,n}$ decreases for $d<n$ and increases for $d\ge n$.  Further, we see that $\mathbb{T}_{n-l-1,n} = \frac{2}{3}\mathbb{T}_{n+l,n}$,
	however computationally  $\mathsf{C}_{n-l-1,n}[f]$ is more efficient than $\mathsf{C}_{n+l,n}[f].$ 
\end{remark}
The following theorem states the error estimate for the truncated Chebyshev series approximation based on the decay estimate \eqref{cjXiang.eq}.
\begin{theorem}\label{thm:errestxiang18}
	Assume the hypotheses of Theorem \ref{thm:decayc_jinV_k} for an integer $k\ge 1$. Let for given integers $d$ and $n$ such that $n\geq k$ and $k\leq d \leq 2n-k-1$, $\mathsf{C}_{d,n}[f]$ be the truncated Chebyshev series of $f$ with approximated coefficients. Then we have
	\begin{equation}
		\|f-\mathsf{C}_{d,n}[f]\|_1 \leq \mathbb{T}_{d,n},
	\end{equation}
	where
	\begin{enumerate}
		\item if $d = n-l$, for some $l = 1,2,\ldots,n-k$, we have 
		\begin{equation}
			\mathbb{T}_{d,n} =\dfrac{4V_k(b-a)^{k+1}}{4^{k}k\pi}\sum_{j=0}^{k}{k \choose j}\dfrac{1}{(n-l+j)(n-l+j-1)\cdots(n-l+j-k+1)}
		\end{equation}
		\item if $d = n+l$, for some $l = 0,1,2,\ldots,n-k-1$, we have 
		\begin{equation}
			\mathbb{T}_{d,n}=\dfrac{6V_k(b-a)^{k+1}}{4^{k}k\pi}\sum_{j=0}^{k}{k \choose j}\dfrac{1}{(n-l+j-1)(n-l+j-2)\cdots(n-l+j-k)}.
		\end{equation}
	\end{enumerate}   
\end{theorem}
\textbf{Proof:} Recall the $L^1$-error estimate for the truncated Chebyshev series expansion of $f$ with approximated coefficients given by \eqref{Est01.eq} and \eqref{E.estimate}, 
\begin{equation}\label{eq:C_d,n1normrestrictedTemp}
	\|f-\mathsf{C}_{d,n}[f]\|_1 \leq \begin{cases}
		2(b-a)\sum\limits_{i=n-l+1}^{\infty}|c_i|, &\quad\text{for } d=n-l,l=1,2,\ldots,n-k,\\
		3(b-a)\sum\limits_{i=n-l}^{\infty}|c_i|, &\quad\text{for } d=n+l,l=0,1,\ldots,n-k-1.
	\end{cases}
\end{equation}
{\bf Case 1:} Now let us take the first case in the above estimate and apply  Theorem \ref{thm:decayc_jinV_k} for $d = n-l,l=1,2,\ldots,n-k$ to get
\begin{align*}
	\sum_{i=n-l+1}^{\infty}|c_i| & \leq \dfrac{2V_k}{\pi}\left(\dfrac{b-a}{4}\right)^{k}\sum_{j=0}^{k}{k \choose j}\sum_{i=n-l+1}^{\infty}\dfrac{1}{(i+j)(i+j-1)\ldots (i+j-k)}\\
	& = \dfrac{2V_k}{\pi}\left(\dfrac{b-a}{4}\right)^{k}\sum_{j=0}^{k}{k \choose j}\sum_{i=n-l+1}^{\infty}\dfrac{1}{k}\left(\dfrac{1}{(i+j-1)\cdots(i+j-k)}-\dfrac{1}{(i+j)\cdots(i+j-k+1)}\right)\\
\end{align*}
Hence using telescopic property of the above series, we have
\begin{align}\label{eq:sumd=n-l}
	\sum_{i=n-l+1}^{\infty}|c_i| & \leq \dfrac{2V_k}{k\pi}\left(\dfrac{b-a}{4}\right)^{k}\sum_{j=0}^{k}{k \choose j}\dfrac{1}{(n-l+j)(n-l+j-1)\cdots(n-l+j-k+1)}.
\end{align}
{\bf Case 2:} Similarly, for $d=n+l,l=0,1,\ldots,n-k-1$, we apply  Theorem \ref{thm:decayc_jinV_k} to get
\begin{equation}\label{eq:sumd=n+l}
	\sum_{i=n-l}^{\infty}|c_i| \leq \dfrac{2V_k}{k\pi}\left(\dfrac{b-a}{4}\right)^{k}\sum_{j=0}^{k}{k \choose j}\dfrac{1}{(n-l+j-1)(n-l+j-2)\cdots(n-l+j-k)}.
\end{equation}
On substitution \eqref{eq:sumd=n-l} and \eqref{eq:sumd=n+l} in \eqref{eq:C_d,n1normrestrictedTemp}, we get 
\begin{equation*}
	\|f-\mathsf{C}_{d,n}[f]\|_1 \leq \begin{cases}
		2(b-a)\dfrac{2V_k}{k\pi}\left(\dfrac{b-a}{4}\right)^{k}\displaystyle\sum_{j=0}^{k}{k \choose j}\dfrac{1}{(n-l+j)(n-l+j-1)\cdots(n-l+j-k+1)},\\ \quad\text{for } d=n-l,l=1,2,\ldots,n-k,\\
		3(b-a)\dfrac{2V_k}{k\pi}\left(\dfrac{b-a}{4}\right)^{k}\displaystyle\sum_{j=0}^{k}{k \choose j}\dfrac{1}{(n-l+j-1)(n-l+j-2)\cdots(n-l+j-k)}, \\
		\quad\text{for } d=n+l,l=0,1,\ldots,n-k-1.
	\end{cases}
\end{equation*} 
Hence the required results.
\section{Numerical Comparison}\label{sec:numericalcomdecayest}
In this section, we numerically illustrate that the improved decay estimate of the Chebyshev coefficients and the $L^1$-error estimate of the truncated Chebyshev series approximation we obtain in section \ref{Cheb.Approx.sec}, are sharper than the earlier ones, obtained in \cite{maj_17a}.
\begin{example}\label{ex:decayest}
	Let us consider the following example
	\begin{equation}\label{eq:decayestex}
		g(t) = \dfrac{|t|}{t+2}, ~~~~~~~ t\in [-1,1].
	\end{equation}
\end{example}
The function $g$ is absolutely continuous and 
$$
g'(t) = \left\{\begin{array}{lcr}
	\dfrac{-2}{(t+2)^2},&~\mbox{ if }& -1\le t < 0\\
	\dfrac{2}{(x+2)^2},&~\mbox{ if }& 0<t\le 1,
\end{array}\right.
$$
which is not continuous. Therefore, we have to take $k=1$. Let us check the other hypothesis of Theorem \ref{thm:coeffbound} and Theorem \ref{thm:decayc_jinV_k}.  Using the weak derivative of $g'$, we can compute $V_1$ in Theorems \ref{thm:errestimate} and  \ref{thm:errestxiang18} as
$$V_1 = 1+ \dfrac{2\pi}{\sqrt{3}} < \infty$$
and the bounded variation of $g'$ is approximately equal to 2.7778, which is taken as the value of $V_1$ in Theorems \ref{thm:decayc_jinV_k} and \ref{thm:errestxiang18}.
\begin{figure}[t]
	\centering
	\includegraphics[height=7cm,width=7.5cm]{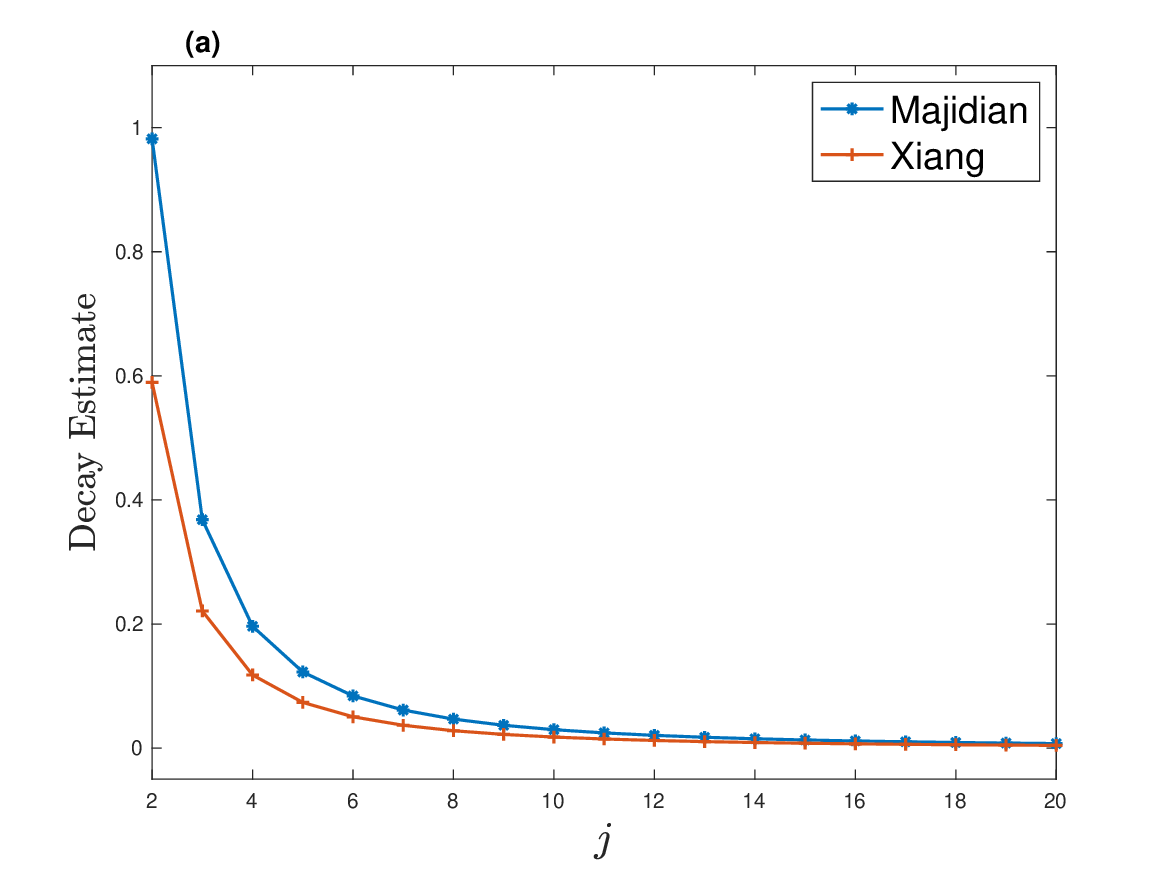}
	\includegraphics[height=7cm,width=7.5cm]{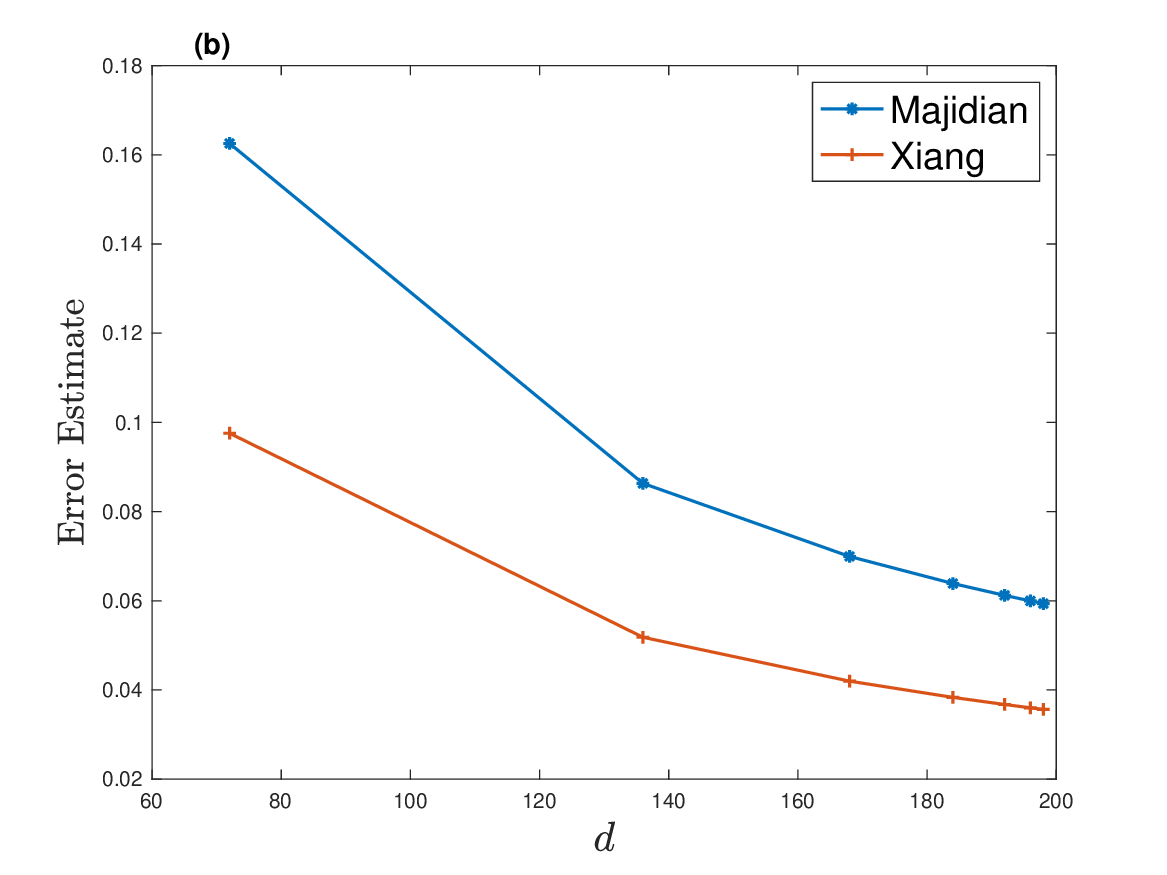}
	\caption {(a) Depicts the comparison between the decay bounds of $|c_j|$, for $j=2,3,\ldots, 30$ derived in Theorem \ref{thm:errestimate} (blue line) and Theorem \ref{thm:decayc_jinV_k} (red line). (b) Depicts the comparison between the error estimates we obtain in Theorem \ref{thm:errestimate}  (blue line) and Theorem \ref{thm:errestxiang18} (red line) for $d=n-l$, where $n=200$ and $l=2^j,j=1,2,\ldots,7$. }
	\label{fig:deacyest}
\end{figure}
The decay estimates (bounds given in \eqref{eq:coeffbound1} and \eqref{cjXiang.eq}) of the Chebyshev series coefficients $c_j$ of $g$, for $j=2,3,\ldots$,  as a function $j$ given in Theorem \ref{thm:errestimate} and Theorem \ref{thm:decayc_jinV_k} are depicted in Figure \ref{fig:deacyest} (a). The error estimates for the truncated Chebyshev series that we obtained in Theorem \ref{thm:errestimate} and Theorem \ref{thm:errestxiang18}, are demonstrated in Figure \ref{fig:deacyest}(b). It can be seen that the improved estimates we derived by using the results of Xiang \cite{xia-18a} are sharper than the earlier ones we obtained using Majidian's \cite{has-17a}. 
\section{Conclusion and Future Directions} $\!\!\!$
\section{Conclusion and Future Directions}
In conclusion, this study has successfully derived generalized decay bounds for Chebyshev series coefficients applicable to functions defined on intervals beyond the standard $[-1,1]$ domain. This generalization holds significance across various applications, particularly in scenarios where the function to be approximated extends beyond the conventional interval boundaries. By leveraging these decay estimates, we have also established truncation error results for a broader class of Chebyshev approximations, considering approximated coefficients of functions $f \in L^2_\omega([a, b])$. This is particularly pertinent as exact computation of Chebyshev coefficients is often impractical in real-world applications, rendering many existing results inapplicable. The optimal error bounds derived in this study are poised to impact several areas of contemporary research significantly. For instance, in the realm of Graph Neural Networks (GNNs), where approximating frequency response is a crucial task, Chebyshev approximation emerges as a natural choice due to its intrinsic properties and computational efficiency. Our findings, applicable to a general class of functions characterized by bounded variation, hold promise for enhancing the accuracy and efficiency of approximation techniques across various domains.


	\bibliography{bibfile}
	\appendix 
	\section{} 

\end{document}